\documentclass{article}
\usepackage{amssymb}
\usepackage{amsthm}
\usepackage{verbatim}
\usepackage{enumerate}
\usepackage{comment}

\newtheorem{thm}{Theorem}[section]
\newtheorem{lem}[thm]{Lemma}

\newtheorem{ex}[thm]{Example}

\newcommand{\sii} {\sigma_{-1}}

\title%[Prime Factors of an Odd Perfect Number]
    {New techniques for bounds on the total number of Prime Factors of an 
     Odd Perfect Number}
\author{Kevin G. Hare
\footnote{Department of Pure Mathematics, University of Waterloo, Waterloo,
    Ontario, Canada,  N2L 3G1, kghare@math.uwaterloo.ca}}

\date{April 12, 2006}

\begin{document}
\maketitle
\thanks{Research of K.G. Hare supported, in part by NSERC of Canada.}

\begin{abstract}
Let $\sigma(n)$ denote the sum of the positive divisors of $n$.
We say that $n$ is perfect if $\sigma(n) = 2 n$.
Currently there are no known odd perfect numbers.
It is known that if an odd perfect number exists, then it must be 
    of the form $N = p^\alpha \prod_{j=1}^k q_j^{2 \beta_j}$, where
    $p, q_1, \cdots, q_k$ are distinct primes and 
    $p \equiv \alpha\equiv 1 \pmod{4}$.
Define the total number of prime factors of $N$ as 
    $\Omega(N) := \alpha + 2 \sum_{j=1}^k \beta_j$.
Sayers showed that $\Omega(N) \geq 29$.
This was later extended by Iannucci and Sorli to show that $\Omega(N) \geq 37$.
This was extended by the author to show that $\Omega(N) \geq 47$.
Using an idea of Carl Pomerance this paper extends these results.
The current new bound is $\Omega(N) \geq 75$.
\end{abstract}

\section{Introduction}

Here and throughout, $n$ is any natural number, and $N$ is a hypothetical 
    odd perfect number.
Let $\sigma(n)$ denote the sum of the positive divisors of $n$.
We say that $n$ is perfect if $\sigma(n) = 2 n$.
It is known that if $\sigma(n) = 2n $ and $n$ is even, then 
    $n = 2^{k-1} (2^k-1)$ where $2^k-1$ is a Mersenne prime.
Currently there are no known odd perfect numbers.
First shown by Euler, it is well known that if 
    an odd perfect number exists, then it must be 
    of the form 
\begin{eqnarray}
N = p^\alpha \prod_{j=1}^k q_j^{2 \beta_j},
\label{eq:odd}
\end{eqnarray} where
    $p, q_1, \cdots, q_k$ are distinct primes and 
    $p \equiv \alpha\equiv 1 \pmod{4}$.

Based on (\ref{eq:odd}) we define the total number of prime factors 
    of an odd perfect number as
\begin{eqnarray}
\Omega(N) := \alpha + 2 \sum_{j=1}^k \beta_j,
\label{eq:Omega}
\end{eqnarray}
and we define the total number of distinct prime factors of $N$ as
    \begin{eqnarray} \omega(N) := 1 + k. \label{eq:omega} \end{eqnarray}
A number of bounds have been derived for $\Omega(N)$.
Cohen showed that $\Omega(N) \geq 23$ \cite{Cohen82}.
Sayers showed that $\Omega(N) \geq 29$ \cite{Sayers86}.
Iannucci and Sorli showed that $\Omega(N) \geq 37$ \cite{IannucciSorli03}.
The author extended this to give $\Omega(N) \geq 47$ \cite{Hare05a}.
This paper extends this result to give

\begin{thm} If $N$ is an odd perfect number, then $\Omega(N) \geq 75$. 
\label{thm:main}
\end{thm}

In proving these results, the methods of \cite{Hare05a} were modified and 
    applied.  
For an introduction and explanation of the algorithm we refer the 
    interested reader there.

To some extent, with the modification of the algorithm given here, 
    the calculation becomes a matter of book-keeping.
In \cite{Hare05a} there was a specific road-block which prevented any
    further calculation. 
This is no longer the case.  
With enough computational power, we could get to any number, although the
    amount of time required appears to be exponential, more than 
    doubling every time we increase the bound by 4 (see Table \ref{tab:times}).
So it is just a matter of how much computer power to dedicate 
    towards the problem.
The choice of 75 in the theorem was motiviated by the fact that 
    it was sufficiently large to demonstrate 
    the improvement in the algorithm while still remaining reasonable 
    with respect to computation time.
\begin{table}
\caption{Timing results for algorithm}
\label{tab:times}
\begin{tabular}{|l|l|l|}
\hline
Result & Time to Prove result & Lines to Prove result\\
\hline
$\Omega(N) \geq 69$ & 
%                      185 hours, 23 minutes 40 seconds \\
                      7 days, 17 hours, 23 minutes, 50 seconds &
%92:52:18.43 + 92:31:21.93  %%% ALTERED FOR 11
%92+1/60*52+1/60^2*18.43 + 92+1/60*31+1/60^2*21.93
%185.3945445  
25,234,392 lines \\  %% ALTERED FOR 11
$\Omega(N) \geq 71$ & 
%                      262 hours, 17 minutes, 58 seconds  \\
                      10 days, 22 hours, 18 minutes, 4 seconds & 
%131:23:12.35 + 130:54:45.66            %% ALTERED FOR 11
%131+1/60*23+1/60^2*12.35 + 130+1/60*54+1/60^2*45.66;  
% 262.2994472
37,752,127 lines \\  %% ALTERED FOR 11
$\Omega(N) \geq 73$ &  
%                     409 hours, 28 minutes, 28 seconds \\
                       17 days, 1 hour, 28 minutes, 36 seconds & 
%148:10:52.41 + 150:33:25.34a + 2:30:02.21 + 1:57:45.35   %% ALTERED FOR 11
%+ 51:51:06.36 + 45:08:04.90 + 4:19:32.38 + 4:57:39.10
% 409.4744445 total
56,352,999 lines \\   %% ALTERED FOR 11
$\Omega(N) \geq 75$ &
                       25 days, 3 hours, 17 minutes, 35 seconds &  %% ALTERED 11
% 603.290811
83,902,264 lines \\    %% ALTERED 11
%Time.out1
%455.276u 566.316s 42:12:24.80 0.6%      0+0k 0+0io 3pf+0w
%Time.out2
%1637.643u 2082.795s 152:23:21.91 0.6%   0+0k 0+0io 0pf+0w
%Time.out3
%906.597u 1191.579s 79:33:00.86 0.7%     0+0k 0+0io 0pf+0w
%Time.out4
%809.237u 1013.459s 74:05:35.76 0.6%     0+0k 0+0io 2pf+0w
%Time.out5
%927.325u 1160.900s 75:26:58.77 0.7%     0+0k 0+0io 0pf+0w
%Time.out6
%1234.393u 1549.627s 90:21:22.59 0.8%    0+0k 0+0io 5pf+0w
%Time.out7
%1172.815u 1461.271s 89:14:42.23 0.8%    0+0k 0+0io 0pf+0w

%42 + 1/60.*12+1/60.^2*24.80 +
%152 + 1/60.*23+1/60.^2*21.91+
%79 + 1/60.*33+1/60.^2*00.86 +
%74 + 1/60.*05+1/60.^2*35.76 +
%75 + 1/60.*26+1/60.^2*58.77 +
%90 + 1/60*21+1/60^2*22.59 +
%89 + 1/60*14+1/60^2*42.23;

\hline
\end{tabular}
\end{table}

%Between 69 and 71, there is an increase by a factor of about 1.414.
%Between 71 and 73, there is an increase by a factor of about 1.561.
%Between 73 and 75, there is an increase by a factor of about 1.473.
%
%Based on these numbers, we could estimate that the amount of time to compute
%    up to 101 would be over 11 years of CPU time.
%100243.3726
    
\section{Definitions and Notation}

For any prime $p$, by $p^a\parallel N$ we mean $p^a |N$ and $p^{a+1}\nmid N$.
By $p^a\nparallel N$ we mean either $p^{a+1}| N$ or $p^a \nmid N$.
We define the function $\sii(n)$ as 
\begin{eqnarray}
\sii(n) := \sum_{d|n}d^{-1} = \frac{\sigma(n)}{n}.
\label{eq:sii}
\end{eqnarray}

A number of simple results concerning $\sii(n)$ are summarized below.
\begin{lem}
\label{lem:sii}
Let $n$ be any natural number.
Then 
\begin{itemize}
\item $\sii(n)$ is a multiplicative function,  
    i.e. if $(n, m) = 1$ then $\sii(n\cdot m) = \sii(n) \sii(m)$,
\item $\sii(n) > 1$ for all $n > 1$,
\item $\sii(n) = 2$ if and only if $n$ is perfect,
\item $\frac{p+1}{p} \leq \sii(p^a) < \sii(p^{a+1}) < \frac{p}{p-1}$
    for all primes $p$ and integers $a \geq 1$.
\end{itemize}
\end{lem}

There are a number of useful results concerning $\omega(N)$, the total number
    of distinct prime factors.
\begin{lem}
\label{lem:exclude}
Let $N$ be an odd perfect number.
Then
\begin{itemize}
\item $\omega(N) \geq 8$ \cite{Chein79, Hagis80}.
\item If $3 \nmid N$ then $\omega(N) \geq 11$ \cite{Hagis83, Kishore83}.
\item If $3 \nmid N$ and $5 \nmid N$ then $\omega(N) \geq 15$ \cite{Norton61}.
\item If $3 \nmid N$, $5 \nmid N$ and $7 \nmid N$ then
    $\omega(N) \geq 27$ \cite{Norton61}.
\item If $3 \nmid N$, $5 \nmid N$, $7 \nmid N$, and $11 \nmid N$, 
    then $\omega(N) \geq 41$ and hence $\Omega(N) \geq 81$.
\end{itemize}
\end{lem}

The last result is probably not the tightest possible, but it is sufficient
    for our purposes.
It follows immediately by noticing that the product $\prod \frac{p_i}{p_i-1}$ 
    for the first 40 primes strictly greater than 11 is less than 2.
So, for our purposes it suffices to assume that $\Omega(N) \leq 73$ and 
    show as a result that $3\nmid N$, $5\nmid N$, $7\nmid N$ and
    $11 \nmid N$, which gives the desired contradiction.

This last result is well known in the literature.  
We use the version given in \cite{CohenSorli03}.  \nocite{Pomerance74}

\begin{lem}[Cohen, Sorli \cite{CohenSorli03}] \label{lem:Carl}
    Let 
    \[N = p_1^{\alpha_1} p_2^{\alpha_2} \cdots p_k^{\alpha_k}
          q_1^{\beta_1} q_2^{\beta_2} \cdots q_r^{\beta_r} \]
    be an odd perfect number, where the $p_i$ and $q_j$ are distinct primes with
    the $p_i$ and $\alpha_i$ known.
    Define \[ S := \sii\left(p_1^{\alpha_1} \cdots p_k^{\alpha_k}\right). \]

    Then we have,  
    \[ \min q_i < \frac{2 + S (r-1)}{2-S}\]
    moreover, if $r \geq 3$ we have the tighter bound 
    \[ \min q_i < \frac{8 + 2 S^2 (r-1)}{4 - S^2}. \]
    Furthermore we have the lower bound 
    \[ \min q_i \geq \max\left(3, \frac{S}{2 - S}\right). \]
\end{lem}

\begin{proof}
Taking $h = \sii$, $\mu = 1$, $h(\lambda) = S$ and $w = r$, this follows 
    directly from equations (2.3), (4.1) and (2.4) of \cite{CohenSorli03}.
\end{proof}

\begin{ex}
Assume that 
    $N = 3^2\cdot 13^4\cdot 30941 \cdot q_1^{\beta_1} \cdot q_2^{\beta_2}$
    is an odd perfect number.
Then by Lemma \ref{lem:Carl} we can assume that $3.5962 < q_1 < 8.1924$.
So in, particular we can assume that $q_1 = 5$ or $7$.
\end{ex}

\section{The Algorithm and Proof of Theorem \ref{thm:main}}%and \ref{thm:equiv}}

Suppose $N = p^{\alpha}\prod q_i^{2 \beta_i}$, as before.
To prove that $\Omega(N) \geq K$, we assume that $\Omega(N) =
    \alpha + \sum 2 \beta_i \leq K-2$ and obtain a contradiction 
    for every combination of $\alpha$ and $\beta_i$.
As mentioned above, it suffices to assume that $\Omega(N) \leq 73$ and 
    show as a result that $3\nmid N$, $5\nmid N$, $7\nmid N$ and
    $11 \nmid N$ to obtain our contradiction.

There are two main modifications to the algorithm in \cite{Hare05a}.
First, as opposed to doing every individual case of $[\alpha, \beta_1, \cdots
    \beta_k]$ where $\Omega(N) = \alpha + 2 \sum \beta_i$, we combine them into 
    one test.
For example, if we want to prove $3 \nmid N$, where $\Omega(N) \leq 57$, 
    we look at the possibilities 
    $3^2\parallel N, 3^4 \parallel N, \cdots, 3^{56}\parallel N$, and recurse.
(Actually, for any power $\beta \geq 46$, we have that $3^\beta\parallel N$ 
    will give
    rise to a contradiction as $\omega(N) \geq 8$.)
This was done because the original method of looking at every possible 
    partition of $N$ caused a large amount of duplicated effort in the
    automated proofs.
By combining these tests together we speed up the calculation, and 
    significantly reduce the storage space requirements.

It should be pointed out that this modification introduced a new means of 
    obtaining a contradiction,
    which in \cite{Hare05a} was taken care of in the choice of 
    partitions.  
This is listed as contradiction (\ref{cont:omega}) below.

The second modification to the algorithm is the use of Lemma \ref{lem:Carl}.
If the algorithm finds itself in a situation where previously it could not 
    continue, because it did not know the factorization of some very 
    large number, we then compute 
    $\frac{2 + S (r-1)}{2-S}$ or $\frac{8 + 2 S^2 (r-1)}{4 - S^2}$
    as approriate, from Lemma \ref{lem:Carl}.
If this upper bound is reasonably small, then we run though all of the 
    possibilities of $q_i$ prime less than this upper bound, 
   (and greater than the lower bound),
    as a means of continuing the calculation.
Here ``small'' was defined as anything less than 100000.
This was sufficient for these calculations.
It should also be noted that the primes are checked in order. 
After proving that a particular prime (say $p$) causes a contradiction,
    it is assumed that it cannot occur as a factor for the next cases
    being checked, within the same sub-branch.
(For example, with an abuse of notation, we assume $p^0 \parallel N$, and
    would arrive at contradiction (\ref{cont:xs=p}) if a factor of $p$ 
    occurs.)

There are five contradictions that we test for.
The first four are from the original algorithm, or the equivalent variation
    needed to combine all tests into one test.  
The last contradiction is commented on above.
\begin{enumerate}
\item Excess of a given prime: \label{cont:xs=p}

    By assuming $p^k \parallel N $ we derive the contradiction that
     $p^{k+1} | N$.
    This is denoted in the output by ``xs=$p$'' where $p$ is the prime in 
    question.

    For example, if we wish to show that $269 \nparallel N$ assuming
        that $3^2 \parallel N$, we first assume that $269 \parallel N$.
    Then we see that $3^3$ and $5$ must divide $N$ (the factors of 
        $\sigma(269)$), which contradicts $3^2 \parallel N$.
    In this case, this would be denoted in the output by ``xs=3''.

\item Excess of the number of primes: \label{cont:xs=prime}

    We have more primes than we are allowed, given the restrictions on 
    $\Omega(N)$ and the fact that only one prime can have an exponent of 1.
    This is denoted in the output by ``xs=prime''.
    Incompletely factored numbers are counted as contributing 
    two primes, even though this may be too low.
    Incompletely factors numbers are known not to be perfect powers.
    Furthermore, incompletely factored numbers are checked to ensure that 
    they are co-prime with each other, as well as other primes within 
    the relevant branch.

    For example, if we wish to show that $26881^{14} \nparallel N$ 
    when $\Omega(N) \leq 19$, we would start with the
    assumption that $26881^{14} \parallel N$.
    This would imply that 
    3, 5, 31, 43, 3368729516337631, 6717545999551, 5601667, 
    8265157321, 18691, 145861 and 1801 must all divide $N$, as they
    are the factors of $\sigma(26881^{14})$.
    We see that the factorization of $N$ that would maximize the number of prime
        factors would be $N = 26881^{14} \cdot p_1 \cdot p_2^2 \cdot p_3^2$.
    As this has at most 3 primes unassigned, and there are 11 unassigned prime 
        factors of $\sigma(26881^{14})$, we get a contradiction.
   
\item Partition cannot be satisfied: \label{cont:part}

    The factors that must divide $N$, along with their powers, cannot satisfy
        the partition.
    In the original algorithm there were a number of different ways that 
        this contradiction could occur.
    Given that in this implementation, there are not specific exponent
        bounds assigned before starting to recurse, this can only 
        occur in one way.
    This is if, of the remaining primes, one of them has to be the special
        prime (with an exponent $a \equiv 1 \pmod{4}$), and 
        of all of the remaining primes, all of them are such that
        $p \not\equiv 1 \pmod{4}$.
    This is denoted in the output by ``exponent bounds exceeded''.

    For example, if we wished to show that $3^{10} \nparallel N$ when
         $\Omega(N) \leq 13$, we would start by assuming that
         $3^{10} \parallel N$.
    This implies that $23$ and $3851$ both must divide $N$, (the 
         two factors of $\sigma(3^{10})$).
    At this point, as $\Omega(N) \leq 13$, we must have one of these
        two primes being the special prime to get either 
        $N = 3^{10} \cdot 23 \cdot 3851^2$ or 
        $N = 3^{10} \cdot 23^2 \cdot 3851$.
    But we notice that $23 \equiv 3851 \equiv 3 \pmod{4}$, hence 
        we see that neither of 23 or 3851 could be the special prime, 
        hence a contradiction. 

\item Excess of $\sii$:  

    A lower bound for $\sii(N)$ using known factors gives
    $\sii(N) > 2$.
    This is denoted in the code by ``S= number'', giving a floating 
        point approximation for a lower bound of $\sii(N)$.
    (The code uses exact rational arithmetic to check the inequality.)

    For example, if we wished to show that $90089 \nparallel N$ we would
        start by assuming that $90089 \parallel N$.
    This implies that $3^2, 5, 7, 11$ and $13$ all divide $N$
        (the factors of $\sigma(90089)$).
    This implies that 
        \begin{eqnarray*}
      \sii(N) &\geq & \frac{\sigma(3^2)}{3^2}\times
                         \frac{\sigma(5)}{5}\times
                         \frac{\sigma(7)}{7}\times
                         \frac{\sigma(11)}{11}\times
                         \frac{\sigma(13)}{13}\times
                         \frac{\sigma(90089)}{90089} \\
                  & \approx& 2.327298560 \\ &>& 2 
        \end{eqnarray*}
        which is the desired contradiction.
    This is denoted in the output as ``S=2.327298560''.

\item $\omega$-bound exceeded. \label{cont:omega}

    Given the choices of primes and exponents made so far, it is not
        possible to satisfy Lemma \ref{lem:exclude}.

    For example, if we were trying to prove that $3^{46}\nparallel N$ where 
        $\Omega(N) \leq 57$ we would start by assuming that 
        $3^{46} \parallel N$.
    We notice that $\omega(N)$ is maximized by writting 
        $N = 3^{46}\cdot p_1 \cdot p_2^2 \cdot p_3^2
                   \cdot p_4^2\cdot p_5^2$, and hence
        $\omega(N) \leq 7$.
    This contradicts $\omega(N) \geq 8$ from Lemma \ref{lem:exclude}.
    This is denoted as ``violate omega bound''.
\end{enumerate}

In Table \ref{tab:ex} an example of running this code is given, and 
    the five possible exceptions.  
These exceptions are indicated by the numbers (1) through (5) on the 
    left hand side of the page.
The use of Lemma \ref{lem:Carl} is indicated with a star, $(*)$, on the 
    left hand side of the page.
Minor formatting has been done to the output, to avoid lines with more 
    than 80 characters.  
(The output of the actual code would have put the line starting with 
    ``It would be nice to know'' all on one line.)

As was done in \cite{Hare05a}, numbers were factored using the 
    ifactor command in MAPLE, with the easy option specified.
If easy factors were not found, then the number was checked in a 
    hints database (which currently contains over 700 completely, or 
    partially factored numbers).

\begin{table}
\caption{Parts of the proof that $\Omega(N) \geq 75$}
\label{tab:ex}
%\verbatiminput{OUTPUT}
\begin{verbatim}
3^6  => 1093
  1093  => 2 547 
    547^2  => 3 163 613 
      163^2  => 3 7 19 67 
        7^2  => 3 19 
          19^2  => 3 127 
(1)         67^2  => 3 7^2 31  xs=7
    ....
    547^18  => 19256021298645399074821884828797791764310604858317 
      19256021298645399074821884828797791764310604858317^2  => 3 
                  3081128010533825683 143739375561423904832226409 
                  279078490021828554466220167808251320188667868911072727 
      ....
      19256021298645399074821884828797791764310604858317^36  => c_1775
        It would be nice to know more factors of
                   sigma(19256021298645399074821884828797791764310604858317^36) 
(*)     By Cohen/Sorli's argument, N has a prime factor between 3 and 18
        Trying each one in turn
        Next prime to try is 5
        5^2  => 31 
        ....
(5)     17^8  => 19 307 1270657  violate omega bound
        Finished Cohen/Sorli's argument
  ....
  1093^2  => 3 398581 
    398581  => 2 17 19 617
    ....
    398581^4  => 5 1866871 2703853428809791 
      5  => 2 3 
      ....
      5^10  => 12207031 
        1866871^2  => 3 19 331 184725139 
(4)       19^2  => 3 127  S=2.001549342
          ....
          19^16  => 3044803 99995282631947 
            331^2  => 3 7 5233  S=2.310779791
            ...
            331^18  => 282349518620419 8080363351283173500917653281271 
(3)           3044803^2  => 3 3090276117871   exponent bounds exceeded
(2)           3044803^4  => 11 631 12382686952067349629581  xs=prime
....
\end{verbatim}
\end{table}

\section{Comments and Acknowledgments}

In \cite{Hare05a} it was said that there were three numbers
    that needed to be factored, without which the algorithm in 
    \cite{Hare05a} could not be continued.
Using the methods of Carl Pomerance, this problem was avoided. 
At this point, the calculation was terminated not by a particular
    obstruction, but instead because 75 was esthetically pleasing.
Further, 75 was sufficiently large to demonstrate 
    the improvement in the algorithm while still remaining reasonable 
    with respect to computation time.
The purpose of this paper was to demonstrate a new technique to 
    extend this bound, and not necessarily to extend this bound to the 
    farthest extent possible. 

I would like to thank Phil Carmody, Christophe Clavier, Don Leclair
    Paul Leyland, Tom Womack, Paul Zimmermann, and the people at 
    {\tt http://mersenneforum.org} who helped provide factorizations
    for the numbers 
    \begin{itemize}
      \item $\sigma(\sigma(547^{18})^{16})$,  a 789-digit number,
      with a factor of: \[ 1520135498523547561326750429418247. \] 
      \item $\sigma(\sigma(3221^{12})^{22})$, a 927-digit number,
      with a factor of: \[ 46973400441039677515399714233826061. \]
    \end{itemize}
    which I had listed in \cite{Hare05a} as obstructions to my 
    calculations.
This information was included in the hints database used by the code.

I am also indebted to Carl Pomerance for bringing a weaker form of 
    Lemma \ref{lem:Carl} to my attention \cite{Pomerance74}, 
    and more importantly how this could be applied to 
    improve the existing bound.
I would also like to thank the unknown referee for the reference
    \cite{CohenSorli03}, which replaced Lemma 2.3, as well as many other 
    useful suggestions and comments.

Also, I would like to acknowledge William Lipp, the creator of 
   www.oddperfect.org, for providing some additional factorizations.
(The website www.oddperfect.org is a website, which when complete, will 
    co-ordinate an attack to improve the lower bound on odd perfect numbers,
    of $N \geq 10^{300}$, 
    given by \cite{BrentCohenTeRiele91}.)

%\bibliographystyle{amsplain}
%\bibliography{paper}

\begin{thebibliography}{10}

\bibitem{BrentCohenTeRiele91}
R.~P. Brent, G.~L. Cohen, and H.~J.~J. te~Riele, \emph{Improved techniques for
  lower bounds for odd perfect numbers}, Math. Comp. \textbf{57} (1991),
  no.~196, 857--868. \MR{1094940}

\bibitem{Chein79}
E.~Z. Chein, \emph{An odd perfect number has at least 8 prime factors}, Ph.D.
  thesis, Pennsylvania State University, 1979.

\bibitem{Cohen82}
Graeme~L. Cohen, \emph{Generalised quasiperfect numbers}, Ph.D. thesis,
  University of New South Wales, 1982.

\bibitem{CohenSorli03}
Graeme~L. Cohen and Ronald~M. Sorli, \emph{On the number of distinct prime
  factors of an odd perfect number}, J. Discrete Algorithms \textbf{1} (2003),
  no.~1, 21--35, Combinatorial algorithms. \MR{2016472}

\bibitem{Hagis80}
Peter Hagis, Jr., \emph{Outline of a proof that every odd perfect number has at
  least eight prime factors}, Math. Comp. \textbf{35} (1980), no.~151,
  1027--1032. \MR{81k:10004}

\bibitem{Hagis83}
\bysame, \emph{Sketch of a proof that an odd perfect number relatively prime to
  {$3$} has at least eleven prime factors}, Math. Comp. \textbf{40} (1983),
  no.~161, 399--404. \MR{85b:11004}

\bibitem{Hare05a}
Kevin~G. Hare, \emph{More on the total number of prime factors of an odd
  perfect number}, Math. Comp. \textbf{74} (2005), no.~250, 1003--1008
  (electronic). \MR{2114661}

\bibitem{IannucciSorli03}
D.~E. Iannucci and M.~Sorli, \emph{On the total number of prime factors of an
  odd perfect number}, Math. Comp. \textbf{72} (2003), no.~244, 2077--2084.

\bibitem{Kishore83}
Masao Kishore, \emph{Odd perfect numbers not divisible by {$3$}. {II}}, Math.
  Comp. \textbf{40} (1983), no.~161, 405--411. \MR{84d:10009}

\bibitem{Norton61}
Karl~K. Norton, \emph{Remarks on the number of factors of an odd perfect
  number}, Acta Arith. \textbf{6} (1960/1961), 365--374. \MR{26 \#4950}

\bibitem{Pomerance74}
Carl Pomerance, \emph{Odd perfect numbers are divisible by at least seven
  distinct primes}, Acta Arith. \textbf{25} (1973/74), 265--300.
  \MR{0340169}

\bibitem{Sayers86}
M.~Sayers, \emph{An improved lower bound for the total number of prime factors
  of an odd perfect number}, Master's thesis, New South Wales Institute of
  Technology, 1986.

\end{thebibliography}

\providecommand{\bysame}{\leavevmode\hbox to3em{\hrulefill}\thinspace}
\providecommand{\MR}{\relax\ifhmode\unskip\space\fi MR }
\providecommand{\MRhref}[2]{%
  \href{http://www.ams.org/mathscinet-getitem?mr=#1}{#2}
}
\providecommand{\href}[2]{#2}

\end{document}